\documentclass[10pt]{article}
\usepackage{amssymb,amsfonts}
\usepackage{amscd}
\usepackage{verbatim}

\newtheorem{theorem}{Theorem}[section]
\newtheorem{lemma}[theorem]{Lemma}
\newtheorem{proposition}[theorem]{Proposition}

\newtheorem{example}[theorem]{Example}
\newtheorem{definition}[theorem]{Definition}

\newtheorem{remark}[theorem]{Remark}

\newenvironment{proof}{\bf Proof. \rm}{$\Box$}

\newcommand{\be}{\begin{equation}}
\newcommand{\ee}{\end{equation}}

\newcommand{\norm}[1]{\Vert #1 \Vert}

\newcommand{\cB}{\mathcal{B}}

\newcommand{\cD}{\mathcal{D}}
\newcommand{\cE}{\mathcal{E}}

\newcommand{\cF}{\mathcal{F}}

\begin{document}
\title{Duality of $W^{\ast}$-correspondences and applications}
\author{Paul S. Muhly\thanks{Supported by grants from the U.S. National Science
Foundation and from the U.S.-Israel Binational Science Foundation.}\\Department of Mathematics\\University of Iowa\\Iowa City, IA 52242\\\texttt{muhly@math.uiowa.edu}
\and Baruch Solel\thanks{Supported by the U.S.-Israel Binational Science Foundation
and by the Fund for the Promotion of Research at the Technion}\\Department of Mathematics\\Technion\\32000 Haifa\\Israel\\\texttt{mabaruch@techunix.technion.ac.il}}

\maketitle
\begin{section}{Introduction}
In \cite{Pi97} Pimsner used
$C^{\ast}$-correspondences to construct and study  a
rich class of $C^{\ast}$-algebras. In our work \cite{MS98}
 we introduced and studied a class of nonselfadjoint
operator algebras, called tensor algebras,
that are constructed from $C^{\ast}$-correspondences and are
subalgebras of Pimsner's $C^{\ast}$-algebras. Later, in \cite{MSp03},
we studied the Hardy algebras (which are the ultraweak closures of
tensor algebras associated with correspondences over  von Neumann
algebras).
Together they form a rich class of operator algebras containing
 a large variety of  algebras (such as analytic
crossed products \cite{MM83, Pet}, noncommutative disc algebras
\cite{Po91}, free semigroup algebras \cite{DP98}, quiver or free
semigroupoid algebras \cite{MS99,KP03} and others). The definition
of a $C^{\ast}$-correspondence and the construction of the tensor
algebra associated with it will be presented in Section 2.

The simplest example of a tensor algebra is the classical disc
algebra $A(\mathbb{D})$ (associated with the correspondence
$\mathbb{C}$ over the algebra $\mathbb{C}$). The study of its
representation theory amounts to the study of contraction
operators (up to unitary equivalence). It turns out that many
ingredients of that theory can
be generalized to the study of the representation theory of
general tensor or Hardy algebras. In this paper we shall deal
mostly with the Hardy algebras. The simplest one is the classical
$H^{\infty}(\mathbb{T})$.

Model theory was originally formulated by Sz-Nagy and Foias and others (see
\cite{NF66}) to study contraction operators on Hilbert space.
In Section 4 we describe how one can define canonical models for
certain representations of the Hardy algebra. This is mainly an
exposition of some of the results in \cite{MSCM}.
It generalizes some of the classical results as well as the
results of Popescu \cite{Po89} for row contractions.
Our goal is to establish a bijective correspondence between the completely
non coisometric representations of the Hardy algebra and
characteristic operator functions.

In the classical theory the characteristic operator function
defining the model is a Schur multiplier; i.e. an $H^{\infty}$ function defined on
$\mathbb{D}$ taking values in $B(\mathcal{E}_1,\mathcal{E}_2)$ (for a pair of Hilbert
spaces $\mathcal{E}_1, \mathcal{E}_2$). When
$H^{\infty}(\mathbb{D})$ is replaced by a general Hardy algebra
(denoted $H^{\infty}(E)$, where $E$ is a $W^{\ast}$-correspondence over
a von Neumann algebra $M$) the ``Schur multipliers" are elements of
another Hardy algebra. That Hardy algebra is associated with a
correspondence that is ``dual" to $E$ in a sense that we shall make
precise in Definition~\ref{dual}.

The concept of the dual of a $W^{\ast}$-correspondence may
be traced back at least to Rieffel's pioneering work on
Morita equivalence for $C^{\ast}$- and $W^{\ast}$-algebras
\cite{R74b}, although the terminology did not appear until
\cite{MSp03}.  Our thinking about the notion was stimulated
by the work of Arveson \cite{Ar89} where he associated a
Hilbert space (i.e. a correspondence over $\mathbb{C}$) to
an endomorphism $\alpha$ of $B(H)$.  Essentially, he showed
that this Hilbert space is dual to the correspondence
$_{\alpha}B(H)$ associated with the endomorphism. Arveson's
construction was generalized by us \cite{MS02} to yield a
corresponence
 $E_{\Theta}$ associated to a
completely positive map $\Theta$ on a general von Neumann algebra $M$.
The correspondence $E_{\Theta}$ is the dual of the correspondence defined by Paschke in \cite{wP73} as a generalization of the $GNS$ construction associated with a state on a $C^{\ast}$-algebra and
studied extensively in \cite{P86}, \cite{AnD90}, \cite{Mi89} and elsewhere.  Nowadays it is known as the \emph{$GNS$-correspondence} associated to $\Theta$ (see \cite{BS00} especially).
 In
fact, it can be shown that the techniques of
\cite{MS02} that yield endomorphic dilations for semigroups of completely
positive maps are, in a sense that may be made precise, dual to the techniques used
 by Bhat and Skiede in \cite{BS00} to achieve their dilation result.  A discussion of the relation between \cite{MS02} and \cite{BS00} appears in \cite{mSp02} and will be developed further in \cite{MSSp02}.
 Duality of correspondences has also proved useful in the analysis of ``curvature" for completely positive maps
\cite{MS03} and in our work on the Hardy algebra of a $W^{\ast}$-correspondence \cite{MSp03}. It is in study of Hardy algebras where duality first appeared in our thinking even before the appearance of \cite{MS02} (although \cite{MSp03} was completed much later).

In Section 3 we define and discuss duality for $W^{\ast}$-correspondences and
present the duality theorem (Theorem~\ref{dual}). In fact, we
shall generalize here some results (and definitions) of
\cite{MSp03} to the context of $W^{\ast}$-correspondences from one
von Neumann algebra $M$ to another, $N$. (In \cite{MSp03} the main results
were proved under the assumption that $M=N$).
This will allow us to present and sketch the proof of
Theorem~\ref{diffduals} using duality to give necessary and sufficient conditions for
two $W^{\ast}$-correspondences to be Morita equivalent (in the sense of
\cite{MS00}).  Full details will be developed in \cite{MSSp02}.

In the next section we introduce the notation and constructions
used throughout the paper.

\end{section}

\begin{section}{Correspondences and operator algebras}
We start by introducing the basic definitions and constructions.
We shall follow Lance \cite{L94} for the general theory of Hilbert
$C^{\ast}$-modules that we shall use. Let $A$ be a
$C^{\ast}$-algebra and $E$ be a right module over $A$ endowed with
a bi-additive map $\langle \cdot,\cdot \rangle : E\times E\rightarrow
A$ (refered to as an $A$-valued inner product) such that, for
$\xi, \eta \in E$ and $a\in A$, $\langle \xi,\eta a\rangle
=\langle \xi,\eta \rangle a$, $\langle \xi,\eta \rangle ^*=\langle
\eta, \xi \rangle$, and $\langle \xi,\xi \rangle \geq 0$, with
$\langle \xi,\xi \rangle =0$ only when $\xi =0$. Also, $E$ is
assumed to be complete in the norm $\norm{\xi}:=
\norm{\langle \xi,\xi \rangle}^{1/2}$. We write $\mathcal{L}(E)$
for the space of continuous, adjointable, $A$-module maps on $E$.
It is known to be a $C^{\ast}$-algebra. If $M$ is a von Neumann
algebra and if $E$ is a Hilbert $C^{\ast}$-module over $M$, then
$E$ is said to be \emph{self-dual} in case every continuous
$M$-module map from $E$ to $M$ is given by an inner product with
an element of $E$. Let $A$ and $B$ be $C^{\ast}$-algebras. A
\emph{$C^{\ast}$-correspondence} from $A$ to $B$ is a Hilbert
$C^{\ast}$-module $E$ over $B$ endowed with a structure of a left
module over $A$ via a nondegenerate *-homomorphism $\varphi : A\rightarrow
\mathcal{L}(E)$.

When dealing with a specific $C^{\ast}$-correspondence, $E$, from a
$C^{\ast}$-algebra $A$ to a $C^{\ast}$-algebra $B$, it will be convenient to
suppress the $\varphi$ in formulas involving the left action and simply write
$a\xi$ or $a\cdot\xi$ for $\varphi(a)\xi$. \ This should cause no confusion in context.

$C^{\ast}$-correspondences should be viewed as generalized $C^{\ast}%
$-homomorphisms. Indeed, the collection of $C^{\ast}$-algebras together with
(isomorphism classes of) $C^{\ast}$-correspondences is a category that
contains (contravariantly) the category of $C^{\ast}$-algebras and (conjugacy
classes of) $C^{\ast}$-homomorphisms. Of course, for this to make sense, one
has to have a notion of composition of correspondences and a precise notion of
isomorphism. The notion of isomorphism is the obvious one: a bijective,
bimodule map that preserves inner products. Composition is ``tensoring'': If
$E$ is a $C^{\ast}$-correspondence from $A$ to $B$ and if
$F$ is a correspondence from $B$ to $C$, then the balanced tensor
product, $E\otimes_{B}F$ is an $A,C$-bimodule that carries
the inner product defined by the formula
\[
\langle\xi_{1}\otimes\eta_{1},\xi_{2}\otimes\eta_{2}\rangle_{E
\otimes_{B}F}:=\langle\eta_{1},\varphi(\langle\xi_{1},\xi_{2}
\rangle_{E})\eta_{2}\rangle_{F}
\]
The Hausdorff completion of this bimodule is again denoted by $E%
\otimes_{B}F$ and is called the \emph{composition} of $E$
and $F$. At the level of correspondences, composition is not
associative. However, if we pass to isomorphism classes, it is. That is, we
only have an isomorphism $(E\otimes F)\otimes G\simeq E\otimes(F\otimes G)$. It is
worthwhile to emphasize here that while it often is safe to ignore the
distinction between correspondences and their isomorphism classes, at times
 the distinction is of critical importance.

In this paper we deal mostly with correspondences over von Neumann
algebras that satisfy some natural additional properties as
indicated in the following definitions.

\begin{definition}
Let $N$ be von Neumann algebra and let $E$ be a Hilbert
$C^{\ast}$-module over $N$. Then $E$ is called a
\emph{Hilbert }$W^{\ast}$\emph{-module} over $N$ in case $E$
is self-dual.
\end{definition}

\begin{remark}\label{topology}
When $E$ is a $W^{\ast}$-module over a $W^{\ast}$-algebra
$N$, then $E$ carries a natural topology making $E$ a dual
space.  This was proved by Paschke in \cite[Proposition
3.8]{wP73}, where he shows that $E$ may be identified with a
weak-* closed subspace of the dual of the projective tensor
product of the complex conjugate of $E$ with the pre-dual of
$N$.  The weak-* topology on $E$ is called the
$\sigma$-\emph{topology}.  It follows easily that
$\mathcal{L}(E)$ carries a natural topology making
$\mathcal{L}(E)$ a dual space \cite[Remark 3.9 and Proposition 3.10]{wP73}.  Thus, since $\mathcal{L}(E)$
is already a $C^{\ast}$-algebra, we see that it is an
abstract $W^{\ast}$-algebra.  For more details about the $\sigma$-topology see \cite{BDH88}.
\end{remark}

\begin{definition}
If $M$ and $N$ are two von Neumann algebras, then an
$M$-$N$-bimodule $E$ is called a
$W^{\ast}$\emph{-correspondence from }$M$\emph{\ to
}$N$\emph{\ }in case $E$ is a self-dual $C^{\ast}
$-correspondence from $M$ to $N$ such that the $\ast
$-homomorphism $\varphi:M\rightarrow\mathcal{L}(E)$ giving
the left module structure on $E$ is normal. (Note: This
makes sense by virtue of the preceeding remark.)

If $M=N$ we shall say that $E$ is a
$W^{\ast}$-correspondence \emph{over $M$}.
\end{definition}

\begin{remark}\label{isomorph}
An isomorphism of a $W^{\ast}$-correspondence $E_1$ from
$M_1$ to $N_1$ and a  $W^{\ast}$-correspondence $E_2$ from
$M_2$ to $N_2$ is a triple $(\sigma,\Psi,
\tau)$ where $\sigma: M_1 \rightarrow M_2$ and
$\tau: N_1 \rightarrow N_2$ are isomorphisms
of von Neumann algebras, $\Psi: E_1 \rightarrow E_2$ is a vector
space isomorphism preserving the $\sigma$-topology and
for ${\xi}$,$\eta \in E_1$ and $a\in M_1, b\in N_1$, we
have
$\Psi(a{\xi}b)=\sigma(a)\Psi({\xi})\tau(b)$ and $\langle \Psi({\xi}),\Psi({\eta})
\rangle = \tau(\langle {\xi},{\eta} \rangle ) $.
When dealing with correspondences over $M$ and over
$N$ (i.e. when $M_i=N_i,\;i=1,2$),
we shall require that $\sigma=\tau$ (unless we say otherwise).
\end{remark}

It is evident that the composition of two $W^{\ast}$-correspondences is a
$C^{\ast}$-correspondence. However, it is not in general a $W^{\ast}$-correspondence.  One must form its ``self-dual completion".  A few words about this may be helpful.  Suppose $Z$ is a Hilbert module over a von Neumann algebra $N$. Then, as Paschke showed in Theorem 3.2 of \cite{wP73}, $Z$ may be embedded in a $W^{\ast}$-module $X$ over $N$ in such a way that $Z$ is dense in $X$ with respect to the $\sigma$-topology.  Further, $X$ is unique up to isomorphism and so may be referred to as \emph{the self-dual completion of $Z$}.  Paschke's proof requires a passage to the Banach space double dual of the von Neumann algebra $N$.  In Proposition 6.10 of \cite{R74b}, Rieffel gives an alternate approach to the notion of the self-dual completion that works as well for correspondences.  That is, he shows that if $M$ and $N$ are von Neumann algebras and if $Z$ is a $C^{\ast}$-correspondence from $M$ to $N$, then there is an essentially unique $W^{\ast}$-correspondence $X$ from $M$ to $N$ that contains $Z$ as a subspace that is dense in the $\sigma$-topology.  Rieffel's proof also uses "duality" techniques, but of the kind that we discuss below.

\begin{definition}\label{tensor product}
Let $M$, $N$ and $P$ be three von Neumann algebras, let $Y$ be a $W^{\ast}$-correspondence from $M$ to $N$ and let $Z$ be a $W^{\ast}$-correspondence from $N$ to $P$.  Then the self-dual completion of the $C^{\ast}$-correspondence from $M$ to $P$, $Y\otimes_{N}Z$, is called \emph{the $W^{\ast}$-tensor product of $Y$ and $Z$ (balanced over $N$)} and is also denoted $Y\otimes_{N}Z$.
\end{definition}

Since we will be considering only $W^{\ast}$-tensor products of $W^{\ast}$-correspondences in this note, there should be no confusion caused from the potential dual use of the notation $Y\otimes_{N}Z$.

Note also, in particular, that a
$W^{\ast}$-correspondence from a von Neumann algebra $N$ to
$\mathbb{C}$ is a Hilbert space $H$ equipped with a (normal)
representation of $N$. If $E$ is a $W^{\ast}$-correspondence from
$M$ to $N$ then the $W^{\ast}$-tensor product, $E\otimes_{N} H$, is a
Hilbert space equipped with a normal representation of $M$. If
$\sigma$ is the representation of $N$ on $H$ we shall also write
$E\otimes_{\sigma} H$ for this tensor product. Observe that
given an operator $X\in \mathcal{L}(E)$ and an operator $S\in
\sigma(N)'$, the map $\xi \otimes h \mapsto X\xi \otimes Sh$
defines a bounded operator on $E\otimes_{\sigma} H$ denoted by
$X\otimes S$.  This is a consequence of Theorem 6.23 in \cite{R74}.

Observe that if $E$ is a $W^{\ast}$-correspondence over a von Neumann algebra
$M$, then each of the tensor powers of $E$, is a $W^{\ast}$-correspondence over $M$ and so, too, is the full Fock space $\mathcal{F}%
(E)$, which is defined to be the direct sum $M\oplus E\oplus E^{\otimes
2}\oplus\cdots$, with its obvious structure as a right Hilbert module over $M$
and left action given by the map $\varphi_{\infty}$, defined by the formula
$\varphi_{\infty}(a):=diag\{a,\varphi(a),\varphi^{(2)}(a),\varphi
^{(3)}(a),\cdots\}$, where for all $n$, $\varphi^{(n)}(a)(\xi_{1}\otimes
\xi_{2}\otimes\cdots\xi_{n})=(\varphi(a)\xi_{1})\otimes\xi_{2}\otimes\cdots
\xi_{n}$, $\xi_{1}\otimes\xi_{2}\otimes\cdots\xi_{n}\in E^{\otimes n}$.

The
\emph{tensor algebra} over $E$, denoted $\mathcal{T}_{+}(E)$, is defined to be
the norm-closed subalgebra of $\mathcal{L}(\mathcal{F}(E))$ generated by
$\varphi_{\infty}(M)$ and the \emph{creation operators} $T_{\xi}$, $\xi\in E$,
defined by the formula $T_{\xi}\eta=\xi\otimes\eta$, $\eta\in\mathcal{F}(E)$.
 We refer the
reader to \cite{MS98} for the basic facts about $\mathcal{T}_{+}(E)$.

\begin{definition}
\label{Hinfty}Given a $W^{\ast}$-correspondence $E$ over the von Neumann
algebra $M$, the ultraweak closure of the tensor algebra of $E$,
$\mathcal{T}_{+}(E)$, in $\mathcal{L}(\mathcal{F}(E))$, will be called the
\emph{Hardy Algebra of }$E$, and will be denoted $H^{\infty}(E)$.
\end{definition}

\begin{example}\label{ex0.1}
If $M=E=\mathbb{C}$ then $\mathcal{F}(E)$ can be identified with
$H^2(\mathbb{T})$. The tensor algebra then is isomorphic to the
disc algebra $A(\mathbb{D})$ and the Hardy algebra is the
classical Hardy algebra $H^{\infty}(\mathbb{T})$.
\end{example}

\begin{example}\label{ex0.2}
If $M=\mathbb{C}$ and $E=\mathbb{C}^n$ then $\mathcal{F}(E)$ can
be identified with the space $l_2(\mathbb{F}_n^+)$ where
$\mathbb{F}_n^+$ is the free semigroup on $n$ generators. The tensor
algebra then is what Popescu refered to as the ``non commutative
disc algebra" $\mathcal{A}_n$ and the Hardy algebra is its
$w^*$-closure. It was studied by Popescu \cite{Po91} and by
Davidson and Pitts who denoted it by $\mathcal{L}_n$ \cite{DP98}.
\end{example}

\begin{example}\label{ex0.3}
Let $M$ be a von Neumann algebra and $\alpha$ be an injective normal
$^{\ast}$-endomorphism on $M$. The correspondence $E$ associated
with $\alpha$ is equal to $M$ as a vector space. The right action
is by multiplication, the $M$-valued inner product is $\langle a,b
\rangle =a^*b$ and the left action is given by $\alpha$; i.e.
$\varphi(a)b=\alpha(a)b$. We write $_{\alpha}M$ for $E$. It is
easy to check that $E^{\otimes n}$ is isomorphic to
$_{\alpha^n}M$. The Hardy algebra in this case will be refered to
as the \emph{non selfadjoint crossed product} of $M$ by $\alpha$
and is related to the algebras studied in \cite{MM83} and \cite{Pet}.
\end{example}

\begin{example}\label{id}
If $\alpha$ is the identity endomorphism of $M$, the correspondence $_{\alpha}M$ is called
\emph{the identity correspondence} over $M$.
\end{example}

\begin{example}\label{cp}
Suppose now that $\Theta$ is a normal, contractive, completely
positive map on a von Neumann algebra $M$. Then we can associate
with it the correspondence $M\otimes_{\Theta}M$ obtained by
defining on the algebraic tenson product $M\otimes M$ the
$M$-valued inner product $\langle a\otimes b,c\otimes d
\rangle=b^*\theta(a^*c)d$ and completing. (The bimodule structure
is by left and right multiplications). This correspondence, first defined by Paschke in \cite{wP73}, was
used by Popa \cite{P86}, Mingo \cite{Mi89},
Anantharam-Delarouche \cite{AnD90} and others to study the map
$\Theta$. (In \cite{BS00} it is refered to as the GNS-module).
 If $\Theta$ is an endomorphism this correspondence is
the one described in Example~\ref{ex0.3}; that is, if $\Theta$ is an endomorphism of $M$, then $M\otimes_{\Theta}M=$ $_{\Theta}M$.
\end{example}

%\begin{example}\label{ex0.4}
%Let $M$ be $D_n$, the diagonal $n\times n$ matrices and $E$ be the
%set of all $n\times n$ matrices $A=(a_{ij})$ with $a_{ij}=0$
%unless $j=i+1$ with the inner product $\langle A,B \rangle =A^*B$
%and the left and right actions given by matrix multiplication.
%Then the Hardy algebra is isomorphic to $T_n$, the $n\times n$
%upper triangular matrices. In fact, a similar argument works to
%show that, for every finite nest of projections $\mathcal{N}$ on a
%Hilbert space $H$, the nest algebra $alg \mathcal{N}$ (i.e. the
%set of all operators on $H$ leaving the ranges of the projections
%in $\mathcal{N}$ invariant ) can be written as $H^{\infty}(E)$ for
%some $W^{\ast}$-correspondence $E$.
%\end{example}

In most respects, the representation theory of $H^{\infty}(E)$ follows the
lines of the representation theory of $\mathcal{T}_{+}(E)$. However, there are
some differences that will be important here. To help illuminate these, we
need to review some of the basic ideas from \cite{MS98, MS99, MS02}.

\begin{definition}
\label{Definition1.12}Let $E$ be a $W^{\ast}$-correspondence over a von
Neumann algebra $M$. Then:

\begin{enumerate}
\item  A \emph{completely contractive covariant representation }of $E$ on a
Hilbert space $H$ is a pair $(T,\sigma)$, where

\begin{enumerate}
\item $\sigma$ is a normal $\ast$-representation of $N$ in $B(H)$.

\item $T$ is a linear, completely contractive map from $E$ to $B(H)$ that is
continuous in the $\sigma$-topology on $E$ (Remark \ref{topology}) and the ultraweak
topology on $B(H)$.

\item $T$ is a bimodule map in the sense that $T(S\xi R)=\sigma(S)T(\xi
)\sigma(R)$, $\xi\in E$, and $S,R\in M$.
\end{enumerate}

\item  A completely contractive covariant representation $(T,\sigma)$ of $E$
in $B(H)$ is called \emph{isometric }in case
\begin{equation}\label{isometric}
T(\xi)^{\ast}T(\eta)=\sigma(\langle\xi,\eta\rangle)
\end{equation}
for all $\xi,\eta\in E$.
\end{enumerate}
\end{definition}

It should be noted that the operator space structure on $E$  to which
Definition \ref{Definition1.12} refers is that which $E$ inherits when viewed
as a subspace of its linking algebra. Also, we shall refer to an isometric,
completely contractive, covariant representation simply as an isometric
covariant representation. There is no problem with doing this because it is
easy to see that if one has a pair $(T,\sigma)$ satisfying all the conditions
of part 1 of Definition \ref{Definition1.12}, except possibly the complete
contractivity assumption, but which is isometric in the sense of equation
(\ref{isometric}), then necessarily $T$ is completely contractive. (See
\cite{MS98}.)

As we showed in \cite[Lemmas 3.4--3.6]{MS98}
and in \cite{MSp03}, if a completely contractive covariant representation,
$(T,\sigma)$, of $E$ in $B(H)$ is given, then it determines a contraction
$\tilde{T}:E\otimes_{\sigma}H\rightarrow H$ defined by the formula $\tilde
{T}(\eta\otimes h):=T(\eta)h$, $\eta\otimes h\in E\otimes_{\sigma}H$. The
operator $\tilde{T}$ intertwines the representation $\sigma$ on
$H$ and the induced representation $\sigma^E:= \varphi (\cdot)
\otimes I_H$ on $E\otimes_{\sigma} H$; i.e.
\begin{equation}\label{covariance}
\tilde{T}(\varphi(\cdot)\otimes I)=\sigma(\cdot)\tilde{T}.
\end{equation}
In fact we have the following lemma from \cite[Lemma 2.16]{MSp03}.

\begin{lemma}
\label{CovRep}The map $(T,\sigma)\rightarrow\tilde{T}$ is a bijection between
all completely contractive covariant representations $(T,\sigma)$ of $E$ on
the Hilbert space $H$ and contractive operators $\tilde{T}:E\otimes_{\sigma
}H\rightarrow H$ that satisfy equation (\ref{covariance}). Given such a
$\tilde{T}$ satisfying this equation, $T$, defined by the formula
$T(\xi)h:=\tilde{T}(\xi\otimes h)$, together with $\sigma$ is a completely
contractive covariant representation of $E$ on $H$. Further, $(T,\sigma)$ is
isometric if and only if $\tilde{T}$ is an isometry.
\end{lemma}

Associated with $(T,\sigma)$ we also have maps $\tilde{T}_n :
E^{\otimes n} \otimes H \rightarrow H$ defined by
$\tilde{T}_n(\xi_1 \otimes \xi_2 \cdots \otimes \xi_n \otimes
h)=T(\xi_1)T(\xi_2)\cdots T(\xi_n)h$.

In the next theorem we give a
a complete description of the representations of the tensor
algebra. The main ingredient in the proof of this theorem is the construction,
to a given completely contractive
covariant representation $(T,\sigma)$ of $E$, of an isometric
representation $(V,\pi)$ that dilates it and is the minimal
isometric dilation. Since we shall later need the notation set in
this construction, we briefly describe it. (For full detailes see
\cite{MS98}).

Given $(T,\sigma)$ on $H$ we set $\Delta = (I-\tilde{T}^{*} \tilde{T})^{1/2}$ (in $
B(E \otimes_{\sigma} H)$), $ \Delta _{*} = (I- \tilde{T}
\tilde{T}^{*})^{1/2} $( in $B(H)$), $ \mathcal D = \overline{ \Delta(E
\otimes_\sigma H)} $ and $\mathcal D_{*} =
\overline{\Delta_{*}(H)}$. Also let $L_{\xi} : H \rightarrow E
\otimes _{\sigma} H$ be the map $L_{\xi} h = \xi \otimes h$ and $
D(\xi) = \Delta \circ L_{\xi} : H \rightarrow E \otimes_{\sigma}
H$. Note that $T(\xi) = \tilde{T} \circ L_{\xi} $.
\\ The
representation space $K$ of $(V, \pi)$ is \[ K = H \oplus
\mathcal D \oplus (E \otimes _{\sigma_{1}} \mathcal D ) \oplus ( E
^{\otimes 2} \otimes _{\sigma_ {1}} \mathcal D ) \oplus ... \]
where $\sigma_1(a)$ is the restriction to $\cD$ of $\varphi(a)
\otimes I_H$. The representation $\pi$ can be defined by
 $\pi = diag(\sigma, \sigma_1, \sigma_2, \ldots )$ where
$\sigma_{k+1}(a)=\varphi_k(a) \otimes I_{\cD}$.
 The map $V : E \rightarrow B(K) $ is defined by
\be\label{vxi} V(\xi) =  \left( \begin{array}{clll} T(\xi) & 0 & 0 &\cdots
\\ D(\xi) & 0 & 0 & \cdots \\ 0&L_{\xi}&0&\\ 0 & 0 & L_{\xi} & \\
& & & \ddots \\ \end{array} \right) \label{V(xi)} \ee where
$L_{\xi}$ here is the obvious map from $E^{ \otimes m} \otimes
\mathcal D$ to $ E^{ \otimes (m+1)} \otimes \mathcal D$.

It turns out that $(V,\pi)$ is a covariant isometric
representation of $E$, it dilates $(T,\sigma)$ in the sense that,
for $\xi\in E$ and $a\in M$, $V(\xi)^*$ and $\pi(a)$ leave $H$ invariant and
their
restrictions to $H$ are equal to $T(\xi)^*$ and $\sigma(a)$ respectively.
 It is also a minimal
dilation (in an obvious sense) and can be shown to be \emph{the
unique} (up to unitary equivalence) minimal isometric dilation of
$(T,\sigma)$.

We now write $\mathcal{Q}_0$ for the space of all vectors perpendicular to
every vector of the form $V(\xi)k$, $\xi \in E$, $k \in K$. Note
 that $\mathcal{Q}_0$ is $\pi(M)$-invariant and, for
every $\xi_1, \xi_2, \ldots \xi_m \in E$ and $k_1, k_2$ in $\mathcal{Q}_0$,
$\langle V(\xi_1)V(\xi_2) \cdots V(\xi_m)k_1,k_2 \rangle =0$. Such
a subspace of $K$ is said to be \emph{wandering}. It is easy to
see that $\mathcal{D}$ is also a wandering subspace. Whenever
$\mathcal{M}\subseteq K$ is a wandering subspace, there is a
unitary operator, denoted $W(\mathcal{M})$, from
$\mathcal{F}(E)\otimes_{\pi | \mathcal{M}} \mathcal{M}$ onto the
$V(E)$-invariant subspace of $K$ generated by $\mathcal{M}$
(denoted $L_{\infty}(\mathcal{M})$). We also note that there
is an isometry $u$ from $\mathcal{Q}_0$ onto $\mathcal{D}_*$ that commutes
with $\pi(a)$ for $a\in M$.(See \cite{MSCM}).

\begin{theorem}
\label{Theorem310MS98}Let $E$ be a $W^{\ast}$-correspondence over a von
Neumann algebra $M$. To every completely contractive covariant representation,
$(T,\sigma)$, of $E$ there is a unique completely contractive representation
$\rho$ of the tensor algebra $\mathcal{T}_{+}(E)$ that satisfies
\[
\rho(T_{\xi})=T(\xi)\;\;\;\xi\in E
\]
and%
\[
\rho(\varphi_{\infty}(a))=\sigma(a)\;\;\;a\in M.
\]
The map $(T,\sigma)\mapsto\rho$ is a bijection between the set of all
completely contractive covariant representations of $E$ and all completely
contractive (algebra) representations of $\mathcal{T}_{+}(E)$ whose
restrictions to $\varphi_{\infty}(M)$ are continuous with respect to the
ultraweak topology on $\mathcal{L}(\mathcal{F}(E))$.
\end{theorem}
\begin{definition}
\label{integratedform}If $(T,\sigma)$ is a completely contractive covariant
representation of a $W^{\ast}$-correspondence $E$ over a von Neumann algebra
$M$, we call the representation $\rho$ of $\mathcal{T}_{+}(E)$ described in
Theorem \ref{Theorem310MS98} the \emph{integrated form} of $(T,\sigma)$ and
write $\rho=\sigma\times T$.
\end{definition}

\begin{remark}
\label{keyproblem}One of the principal difficulties one faces in dealing with
$\mathcal{T}_{+}(E)$ and $H^{\infty}(E)$ is to decide when the integrated
form, $\sigma\times T$, of a completely contractive covariant representation
$(T,\sigma)$ extends from $\mathcal{T}_{+}(E)$ to $H^{\infty}(E)$. This
problem arises already in the simplest situation, vis. when $M=\mathbb{C}=E$.
In this setting, $T$ is given by a single contraction operator on a Hilbert
space, $\mathcal{T}_{+}(E)$ ``is'' the disc algebra and $H^{\infty}(E)$ ``is''
the space of bounded analytic functions on the disc. The representation
$\sigma\times T$ extends from the disc algebra to $H^{\infty}(E)$ precisely
when there is no singular part to the spectral measure of the minimal unitary
dilation of $T$. We are not aware of a comparable result in our
general context but we have some sufficient conditions. One of
them is given in the following lemma. It is not necessary in
general.
\end{remark}

\begin{lemma}\label{contraction} \cite[Corollary 2.14]{MSp03}
If $\norm{\tilde{T}}<1$ then $\sigma \times T$ extends to a
$\sigma$-weakly continuous representation of $H^{\infty}(E)$.
\end{lemma}

\end{section}

\begin{section}{Duality of $W^{\ast}$-correspondences}

In this section we discuss the concept of duality for
$W^{\ast}$-correspondences. As we noted in the introduction, the concept arises implicitly in \cite{R74}. It was used in
\cite{Ar89}, again implicitly, to construct the product system associated to an
$E_0$-semigroup on $B(H)$ and a bit more explicitly in \cite{MS02} where the construction
was extended to $E_0$-semigroups on a general von Neumann algebra.

%%%%%%%%%%%%%%%%%%%%%%%%%%%%%%%%%%%%%%%%%%%%%%%%%%%%%%%%%%%%%%%%%%%%%%%%%%%
%% I don't think we need to belabor the priority issue in any more %%detail here.  So I have commented out the rest of the paragraph.
%%and was studied, explicitly, in \cite{MSp03}. It was proved useful
%%to the study of semigroups of completely positive maps (in
%%\cite{MS02}), in the study of the curvature invariant
%%(\cite{MS03}), in interpolation problems (\cite{MSp03}) and, as we
%%shall describe in the last section, to study canonical models of
%%representations.
%%%%%%%%%%%%%%%%%%%%%%%%%%%%%%%%%%%%%%%%%%%%%%%%%%%%%%%%%%%%%%%%%%%%%%%%%%%

\begin{definition}
Let $E$ be a $W^{\ast}$-correspondence from $M$ to $N$. Let
$\sigma:M\rightarrow B(H)$ and $\tau:N\rightarrow B(K)$ be
 normal representations of the von Neumann
algebras $M$ and $N$. Then the $\tau$-$\sigma$\emph{-dual} of $E$, denoted
$E^{\tau,\sigma}$, is defined to
be
\[
\{\eta\in B(H,E\otimes_{\tau}K)\mid \eta\sigma(a)=(\varphi(a)\otimes
I)\eta,\;a \in M \}.
\]
\end{definition}

An important feature of the dual
 $E^{\tau,\sigma}$ is
that it is a $W^{\ast}$-correspondence - over $\sigma(M)^{\prime}$ - as the following proposition shows. (cf. \cite[Theorem 6.5]{R74b}).

\begin{proposition}\label{corres}
With respect to the actions of $\sigma(M)^{\prime}$ and $\tau(N)^{\prime}$ and the $\sigma
(M)^{\prime}$-valued inner product defined as follows, $E^{\tau,\sigma}$ becomes a
$W^{\ast}$-correspondence from $\tau(N)^{\prime}$ to $\sigma(M)^{\prime}$: For
$Y\in\sigma(M)^{\prime}$, $X\in \tau(N)^{\prime}$
 and $T\in E^{\tau,\sigma}$, $X\cdot T\cdot Y:=(I\otimes X)TY$,
and for $T,S\in E^{\tau,\sigma}$, $\langle T,S\rangle_{\sigma(M)^{\prime}}%
:=T^{\ast}S$.
\end{proposition}

\begin{remark}\label{evaluation}
When $M=N$ (i.e. when $E$ is a $W^{\ast}$-correspondence over $M$) and
when $\tau=\sigma$, we write $E^{\sigma}$ in place of
$E^{\sigma,\sigma}$. The importance of this space for us lies in the fact that it is
closely related to the representations of $E$. In fact,
 the operators in $E^{\sigma}$
whose norms do not exceed $1$ are precisely the adjoints of the operators of the
form $\tilde{T}$ for a covariant pair $(T,\sigma)$. In particular,
every $\eta$ in the \emph{open} unit ball of $E^{\sigma}$
(written $\mathbb{D}(E^{\sigma})$) gives
rise to a covariant pair $(T,\sigma)$ (with $\eta=\tilde{T}^*$)
such that $\sigma \times T$ is a representation of
$H^{\infty}(E)$. Given $X\in H^{\infty}(E)$ we can apply this
representation (associated to $\eta$) to it. The resulting
operator in $B(H)$ will be denoted by $X(\eta^*)$.
In this way, we view every element in the Hardy algebra as a
$B(H)$-valued function on the unit ball $\mathbb{D}(E^{\sigma})$.
This point of view is exploited in \cite{MSp03} to deal with
interpolation problems.
\end{remark}

\begin{example}\label{ex0.12}
Suppose $M=E=\mathbb{C}$ and suppose $\sigma$ is the representation of
$\mathbb{C}$ on some Hilbert space $H$ given by scalar multiplication. Then it is easy to check
that $E^{\sigma}$ is isomorphic to $B(H)$. Fix an $X\in
H^{\infty}(E)$. As we mentioned above, this Hardy algebra is the
classical space $H^{\infty}(\mathbb{T})$ and we can identify $X$ with a
function $f\in H^{\infty}(\mathbb{T})$. Given $S\in E^{\sigma}=B(H)$,
it is not hard to check that $X(S^*)$, as defined above, is the operator that arises through the
$H^{\infty}$-functional calculus, $f(S^*)$.
\end{example}

\begin{example}\label{ex0.7}
If $\Theta$ is a contractive, normal, completely positive map on a
von Neumann algebra $M$ and if $E=M\otimes_{\Theta} M$ (see
Example~\ref{cp} ) then, for every faithful representation
$\sigma$ of $M$ on $H$, the $\sigma$-dual is the space of all
bounded operators mapping $H$ into the Stinespring space $K$
associated with $\Theta$ as a map from $M$ to $B(H)$ that
intertwine the representation $\sigma$ (on $H$) and the
Stinespring representation $\pi$ (on $K$). This correspondence has
proved very useful in the study of completely positive maps. (See
\cite{MS02} and \cite{MS03}). If $M=B(H)$ this is a Hilbert space
and was studied by Arveson \cite{Ar89}.
Note also that, if $\Theta$ is an endomorphism, then this dual
correspondence is the space of all operators on $H$ intertwining
$\sigma$ and $\sigma \circ \Theta$.
\end{example}

We now return to discuss the general case (where $M$ is not
necessarily equal to $N$).
The term ``dual" that we use is justified by the following theorem, which is proved as Theorem 3.6 in \cite{MSp03} under the assumption that $M=N$.  See also \cite{mSp03}. It really is contained in Proposition 6.10 of \cite{R74b}, although the proof given in \cite{MSp03} and outlined below is more elementary and transparent. The full result along with applications to Morita theory (in particular Theorem \ref{diffduals}) will be proved in \cite{MSSp02}.

\begin{theorem}\label{dual} Let $E$ be a $W^{\ast}$-correspondence from $M$ to $N$ and let
$\sigma$, $\tau$
be  faithful, normal representations of $M$ (on $H$) and $N$ (on $K$)
respectively. If we write
$\iota_1$ for the identity representation of $\sigma(M)'$ (on $H$)
and $\iota_2$ for the identity representation of $\tau(N)'$ (on
$K$),
then one can form the $\iota_1$-$\iota_2$-dual of $E^{\tau,\sigma}$ and we have
$$ (E^{\tau,\sigma})^{\iota_1,\iota_2} \cong E .$$

\end{theorem}

The isomorphism in the theorem is $(\sigma,\Psi,\tau)$ from $E$ onto
$ (E^{\tau,\sigma})^{\iota_1,\iota_2}
$ where $\Psi$ is defined by the equation
$$ \Psi(\xi)^*(\eta \otimes h)= L_{\xi}^*\eta(h), $$
with $\xi \in E$, $\eta \in E^{\tau,\sigma}$ and $h \in H$. Note
that $\eta$ is a map from $H$ to $E\otimes K$ so that
$L_{\xi}^*\eta(h)$ lies in $K$ and the equation above defines a
map from $E^{\tau,\sigma}\otimes H$ to $K$ whose adjoint can be
shown to have the intertwining property required from an element
of $ (E^{\tau,\sigma})^{\iota_1,\iota_2}$. In order to prove that
$\Psi$ is onto one uses the following lemma. It was proved in
\cite[Lemma 3.5]{MSp03} for the case $\sigma=\tau$. The proof in the general
case requires only a minor adjustment.

\begin{lemma}\label{span}
When $E$ is as above and $\sigma$ and $\tau$ are faithful
representations of $M$ and $N$ respectively, we have
$$ \bigvee \{ X(H) : \;X\in E^{\tau,\sigma} \} =E\otimes_{\tau} K
.$$
\end{lemma}

The following two lemmas show that the operation of taking duals
``behaves nicely" with respect to direct sums and tensor products.

\begin{lemma}\label{directsum}
Given $W^{\ast}$-correspondences $E_1$ and $E_2$ from $M$ to $N$ and
faithful representations $\sigma$ (of $M$ on $H$) and $\tau$ (of $N$ on $K$)
 we have
 $(E_1 \oplus E_2)^{\tau,\sigma}\cong E_1^{\tau,\sigma} \oplus
E_2^{\tau,\sigma}.$
%\item[(ii)] $(E_1 \otimes E_2)^{\sigma}\cong E_2^{\sigma} \otimes
%E_1^{\sigma}.$
%\item[(iii)] $\mathcal{F}(E)^{\sigma}\cong
%\mathcal{F}(E^{\sigma}).$
%\item[(iv)] The map $\eta \otimes h \mapsto \eta(h)$ induces a
%unitary operator from $E^{\sigma}\otimes_{\iota}H$ onto
%$E\otimes_{\sigma}H$.
%\item[(v)] Applying item (iv) above to $\mathcal{F}(E)$ in place of
%$E$, we get a unitary operator $U$ from $\mathcal{F}(E^{\sigma})\otimes H$
%onto $\mathcal{F}(E)\otimes H$.
%\end{enumerate}
\end{lemma}

\begin{lemma}\label{tensorprod}
Let $E$ be a $W^{\ast}$-correspondence from $M$ to $N$ and $F$ be
a $W^{\ast}$-correspondence from $N$ to $Q$. Let $\sigma$, $\tau$
and $\pi$ be normal faithful representations of $M$, $N$ and $Q$
respectively. Then the map $X\otimes Y \mapsto (I_E \otimes X)Y$
(for $X\in F^{\pi,\tau}$ and $Y\in E^{\tau,\sigma}$)
defines an isomorphism
$$F^{\pi,\tau} \otimes E^{\tau,\sigma} \cong (E\otimes
F)^{\pi,\sigma}.$$
\end{lemma}

Given a $W^{\ast}$-correspondence $Z$ over $M$ and two faithful
representations $\sigma$ and $\tau$ (of $M$), the dual
correspondences $Z^{\sigma}$ and $Z^{\tau}$ are, in general, non
isomorphic but, as we shall show below they are Morita equivalent.
We will also show that the converse holds. For this, we first
recall the definition of Morita equivalence for
$W^{\ast}$-correspondences.

In \cite{MS00} we discussed Morita equivalence for
$C^{\ast}$-correspondences. For $W^{\ast}$-correspondences the
concept is defined similarly (with minor changes). Recall first
that an $M$-$N$ \emph{equivalence bimodule} is an $M$-$N$ bimodule
$X$ that is endowed with $M$- and $N$-valued inner products,
$_M\langle \cdot,\cdot \rangle$ and $\langle \cdot, \cdot
\rangle_N$, making $X$ a full and selfdual (right) Hilbert
$W^{\ast}$-module over $N$ and a full and selfdual (left) Hilbert
$W^{\ast}$-module over $M$ such that
$_M\langle \xi,\eta \rangle \zeta = \xi \langle \eta,\zeta
\rangle_N $ for $\xi,\eta$ and $\zeta$ in $X$.
By definition, the von Neumann algebras $M$ and $N$ are strongly
Morita equivalent in case there is an $M$-$N$ equivalence
bimodule \cite{R74}.

\begin{definition}\label{morita} (cf. \cite{MS00})
 $W^{\ast}$-correspondences $E$ and $F$, over $M$ and $N$
 respectively, are said to be (strongly) Morita equivalent if there is
 an $M$-$N$ equivalence bimodule $X$ such that $E\otimes_M X \cong
 X\otimes_N F$ (where the isomorphism here is a triple
 $(id,\Psi,id)$ for some map $\Psi$).
 \end{definition}

In the statement of the following theorem, the isomorphisms are in
the sense indicated at the end of Remark~\ref{isomorph}.

\begin{theorem}\label{diffduals}
Let $E$ and $F$ be $W^{\ast}$-correspondences over $\sigma$-finite
von Nuemann algebras $M$ and $N$
respectively. Then the following conditions are equivalent.
\begin{itemize}
\item[(1)] There is a  $W^{\ast}$-correspondence $Y$ (over some von
Neumann algebra $Q$) and two faithful representations $\pi_1$ and
$\pi_2$ of $Q$ such that $E\cong Y^{\pi_1}$ and $F\cong
Y^{\pi_2}$.
\item[(2)] The $W^{\ast}$-correspondences $E$ and $F$ are strongly Morita equivalent.
\item[(3)] There are faithful representations $\sigma$ and $\tau$,
of $M$ and $N$ respectively, such that $E^{\sigma}\cong F^{\tau}$.
\end{itemize}
\end{theorem}
\begin{proof} We shall only sketch the proof. For a detailed proof
one needs to be careful about the maps involved in each of the isomorphisms
below (see Remark~\ref{isomorph}). To prove (3) implies
(1) , write $\psi$ for the isomorphism from $\sigma(M)'$ onto
$\tau(N)'$ (implied by the assumption that $E^{\sigma}\cong
F^{\tau}$), write $\iota_1$ and $\iota_2$ for the identity
representations of $\sigma(M)'$ and $\tau(N)'$ respectively and
set $Y=E^{\sigma}$. Then it follows from duality that $E\cong
Y^{\iota_1}$ and $F\cong Y^{\iota_2 \psi}$.

To prove (1) implies (2) let $Z$ be the identity correspondence
$Q$ and set $X=Z^{\pi_1,\pi_2}$. Then $E \otimes X \cong Y^{\pi_1}
\otimes Z^{\pi_1,\pi_2} \cong (Q\otimes Y)^{\pi_1,\pi_2} \cong (Y
\otimes Q)^{\pi_1,\pi_2} \cong Z^{\pi_1,\pi_2} \otimes Y^{\pi_2}
\cong X \otimes F$.

For the last part, assume $E$ and $F$ are Morita equivalent and
$X$ is the equivalence bimodule implementing the equivalence.
Assume $N\subseteq B(K)$ and write $H$ for $X\otimes_N K$. Write $\sigma$
for the identity representation of $N'$ (on $K$) and $\tau$ for
the representation of $N'$ on $H$ defined by $\tau(a)=I_X \otimes
a$. Note that $\tau(N')'=\mathcal{L}(X)\otimes I_K$ and
this algebra is isomorphic to $M$ (since $M\cong \mathcal{L}(X)$
for an equivalence bimodule $X$). Write $\psi$ for this
isomorphism (from $M$ to $\tau(N')'$). Let $\iota_1$ be the identity representation of
$\tau(N')'$ on $H$ and $\iota_2$ be the identity representation of
$N$ on $K$.  Write $Z$ for the identity correspondence of $N'$. Given $x\in
X$, define $S(x)$ to be the map from $K$ to $Z\otimes_{\tau} H$ given
by $S(x)(k)=I\otimes_{\tau}(x \otimes k)$. It is easy to check that
$S(x)$ lies in $Z^{\tau,\sigma}$. In fact, the triple $(\psi, S,id)$
is an isomorphism of $X$ and $Z^{\tau,\sigma}$.
 It then follows from duality that $F^{\iota_2}\cong
F^{\iota_2}\otimes N^{\prime}\cong F^{\iota_2} \otimes
(Z^{\tau,\sigma})^{\iota_2,\iota_1} \cong (Z^{\tau,\sigma}\otimes
F)^{\iota_2,\iota_1} \cong (X \otimes F)^{\iota_2,\iota_1
\psi} \cong (E\otimes X)^{\iota_2, \iota_1 \psi} \cong
(E\otimes_{\psi} Z^{\tau,\sigma})^{\iota_2,\iota_1 \psi} \cong
(Z^{\tau,\sigma})^{\iota_2,\iota_1} \otimes E^{\iota_1 \psi} \cong
N^{\prime} \otimes_{\tau} E^{\iota_1 \psi} \cong E^{\iota_1
\psi}.$

\end{proof}

\end{section}

\begin{section}{Applications of duality: Commutants and Can\-onical models}

In this section, $E$ will be a $W^{\ast}$-correspondence over a
von Neumann algebra $M$.
As was mentioned in the introduction, there are several
applications for the dual correspondences. Here we concentrate on
using the Hardy algebra associated with a dual correspondence in
order to study the Hardy algebra associated with the original
correspondence $E$. The first result along these lines is the
identification of the commutant of $H^{\infty}(E)$ (given in some
induced representation). The other is the development of canonical
models to study representations of $H^{\infty}(E)$. As we shall
see, the characteristic functions can be identified with elements
of $H^{\infty}(E^{\tau})$ for some dual correspondence $E^{\tau}$.

Although $H^{\infty}(E)$ was defined as a subalgebra of
$\mathcal{L}(\mathcal{F}(E))$ it is often useful to consider a
(faithful) representation of it on a Hilbert space. Given a
faithful, normal, representation $\sigma$ of $M$ on $H$ we can
``induce" it to a representation of the Hardy algebra. To do this,
we form the Hilbert space $\mathcal{F}(E)\otimes_{\sigma} H$ and
write
$$ Ind(\sigma)(X)=X\otimes I,\;\; X\in H^{\infty}(E).$$
(Note that $Ind(\sigma)(X)$ is a well defined bounded operator on $\mathcal{F}(E)\otimes_{\sigma} H$ for every $X$ in
$\mathcal{L}(\mathcal{F}(E))$. Such representations were studied
by M. Rieffel in \cite{R74b}). $Ind(\sigma)$ is a faithful
representation and is a homeomorphism with respect to the
$\sigma$-weak topologies.
Similarly one defines $Ind(\iota)$ (where $\iota$ is the identity
representation of $\sigma(M)'$ on $H$) , a representation of
$H^{\infty}(E^{\sigma})$. The following theorem shows that,
roughly speaking, the algebras $H^{\infty}(E)$ and
$H^{\infty}(E^{\sigma})$ are the commutant of each other.
For the proof, see \cite[Theorem 3.9]{MSp03}.

\begin{theorem}\label{commutant} \cite{MSp03}
Let $E$ be a $W^{\ast}$-correspondence over $M$ and $\sigma$ be a
faithful normal representation of $M$ on $H$.
Then there exists a unitary operator $U: \mathcal{F}(E^{\sigma})\otimes H
\rightarrow \mathcal{F}(E)\otimes H$ such that
$$
U^*(Ind(\iota)(H^{\infty}(E^{\sigma})))U=(Ind(\sigma)(H^{\infty}(E)))'
$$
and, consequently,
$$(Ind(\sigma)(H^{\infty}(E)))''=Ind(\sigma)(H^{\infty}(E)). $$
\end{theorem}

\begin{example}\label{quiver}
Given an $n\times n$ matrix $C$ in $M_n(\mathbb{Z}_+)$. One can
associate with it a $W^{\ast}$-correspondence $E(C)$ over the
algebra $D_n$ of all diagonal $n\times n$ matrices. (See
\cite{MS99} for details). The tensor algebra associated with
$E(C)$ is called the \emph{quiver algebra} or the \emph{path
algebra} (associated with the directed graph defined by $C$). It
is a subalgebra of the Cuntz-Krieger $C^{\ast}$-algebra $O_C$.
If $\sigma$ is the identity representation of $D_n$ (on
$\mathbb{C}^n$), then $E(C)^{\sigma}=E(C^t)$ (where $C^t$ is the
transpose matrix). Thus, Theorem~\ref{commutant} gives another
proof of \cite[Proposition 5.4]{MS99} and of \cite[Theorem
5.8]{KP03} . ( Note that a semigroupoid algebra of \cite{KP03} is the
image, under an induced representation, of a quiver algebra.)
\end{example}

Another way in which the Hardy algebra of the dual correspondence
plays a role in studying $H^{\infty}(E)$ is through canonical
models. Here, roughly, the elements of the Hardy algebra of the
dual play the role that ``Schur multipliers" play in the classical theory.

The canonical models are used to study certain representations of
the Hardy algebra.
The representations of $H^{\infty}(E)$ that we shall study here
are the completely noncoisometric (abbreviated c.n.c ) and the
$C_{.0}$-representations. We first recall the definitions. (See
\cite{MSp03} for more details).

 Recall that, given a covariant
representation $(T,\sigma)$ of $E$ on $H$, it has a unique minimal
isometric dilation $(V,\pi)$ on $K$. It was described in Section 2
and we use here the notation introduced there.

\begin{definition}\label{repns}
\begin{itemize}
\item[(1)] A completely contractive covariant representation
$(T,\sigma)$ of $E$ is called a $C_{.0}$-representation if, for
every $h\in H$, $\norm{\tilde{T}_n^* h} \rightarrow 0$.
Equivalently, if $L_{\infty}(\mathcal{Q}_0)=K$.
\item[(2)] A completely contractive covariant representation
$(T,\sigma)$ of $E$ is said to be completely noncoisometric if
$L_{\infty}(\mathcal{D}) \vee L_{\infty}(\mathcal{Q}_0) =K$.
\end{itemize}
\end{definition}

Clearly, every $C_{.0}$-representation is
completely non coisometric.

\begin{theorem}\label{extend} \cite{MSp03}
If $(T,\sigma)$ is a c.n.c representation then $T\times \sigma$
extends to a $w^*$-continuous representation of $H^{\infty}(E)$.
\end{theorem}

We now define characteristic operator functions in this context. It
generalizes the classical case \cite{NF66} and the case studied
by Popescu \cite{Po89}. Here we shall present the constructions
and the main results. Full details will appear in \cite{MSCM}.

\begin{definition}\label{charopfun}
Given a von Neumann algebra $M$ and a W*-correspondence $E$ over
$M$,
A \emph{characteristic operator function} is a tuple
$(\Theta,\cE_1,\cE_2,\tau_1,\tau_2)$ such that
\begin{enumerate}
\item[(i)] For $i=1,2$, $\cE_i$ is a Hilbert space and $\tau_i$ is
a representation of $M$ on $\cE_i$.
\item[(ii)] $\Theta : \mathcal{F}(E)\otimes_{\tau_1} \mathcal{E}_1
\rightarrow \mathcal{F}(E)\otimes_{\tau_2} \mathcal{E}_2$ is a
contraction satisfying
\begin{equation}\label{comm1}
(\varphi_{\infty}(a)\otimes
I_{\mathcal{E}_2})\Theta =\Theta (\varphi_{\infty}(a)\otimes
I_{\mathcal{E}_1}),\;\; a\in M
\end{equation}
and
\begin{equation}\label{comm2}
(T_{\xi}\otimes I_{\mathcal{E}_2})\Theta =\Theta
(T_{\xi}\otimes I_{\mathcal{E}_1}).
\end{equation}
\item[(iii)] There is no non zero vector $x \in \cE_1$ such that
$x=P_{\cE_1}\Theta^*P_{\cE_2}\Theta x$. (We say that $\Theta$ is \emph{purely
contractive}).
\item[(iv)] $\overline{\Delta_{\Theta}(\cF(E)\otimes
\cE_1)}=\overline{\Delta_{\Theta}((\cF(E)\otimes \cE_1)\ominus \cE_1)}$
(where $\Delta_{\Theta}=(I-\Theta^*\Theta)^{1/2}$).

\end{enumerate}

If, in addition, $\Theta$ is an isometry then it
will be called an \emph{inner characteristic operator function}.
(In this case, (iv) holds automatically).
\end{definition}

We shall often refer to $\Theta$ as the characteristic operator
function (when $\mathcal{E}_i$ and $\tau_i$ are assumed to be
known).

Note that it follows from (\ref{comm1}) and (\ref{comm2}) that, if
we write $\mathcal{E}$ for $\mathcal{E}_1\oplus \mathcal{E}_2$ and
let $\tau$ be the representation $\tau_1\oplus \tau_2$, then the
matrix $\left( \begin{array}{cc} 0&0 \\ \Theta &0 \end{array}
\right)$ (viewed as an element of $\mathcal{F}(E)\otimes
\mathcal{E}$)
lies in the commutant of $Ind(\tau)(H^{\infty}(E))$. If
$\tau$ is faithful we will be able to use Theorem~\ref{commutant}
to conclude that the matrix defines an element of $H^{\infty}(E^{\tau})$.
If $\tau$ is not faithful, one can form $\tau'=\tau \oplus \tau_0$
that is a faithful representation (on a larger space) and
consider the $3\times 3$ matrix with $\Theta$ in the $2,1$ entry
instead of the $2\times 2$ matrix above. Since this is just a
technical point, we shall ignore it here and use Theorem~\ref{commutant}
to define $\widehat{\Theta} \in
 H^{\infty}(E^{\tau})$ by
\[
Ind(\iota)(\widehat{\Theta})  = U\left( \begin{array}{cc} 0&0 \\ \Theta &0 \end{array}
\right)U^*.
\]
The left hand side will also be written $\hat{\Theta}\otimes
I_{\mathcal{E}}$. Note that, if $\Theta$ is inner,
$\widehat{\Theta}^*\widehat{\Theta}$ is the projection onto
$\mathcal{E}_1$.
%in fact, given
%$\mathcal{E}_1,\mathcal{E}_2,\tau_1$ and $\tau_2$,
% the map $\Theta \mapsto \widehat{\Theta}$ is
%a bijection between all contractions satisfying condition (ii) of
%Definition~\ref{charopfun} and all
%contractions $\Psi\in H^{\infty}(E^{\tau})$ satisfying
%$\Psi=q_2\Psi q_1$ (where $q_i$ is the projection onto
%$\mathcal{E}_i$, $i=1,2$). Inner characteristic operator functions
%are mapped, under this bijection onto these $\Psi$'s that are
%partial isometries with initial space $\mathcal{E}_1$. We write
%$\Psi \mapsto \check{\Psi}$ for the inverse map.

Next we construct, for a fixed characteristic operator function $\Theta$
 a covariant representation associated to it.

For this, we write $\Delta_{\Theta}=(I_{\cF(E)\otimes \cE_1}-\Theta^*\Theta)^{1/2} \in
B(\cF(E)\otimes \cE_1) $ and set
$$K(\Theta)=(\cF(E) \otimes \cE_2)\oplus \overline{\Delta_{\Theta}(\cF(E)\otimes
\cE_1)} \subseteq \cF(E)\otimes \cE$$ and
$$H(\Theta)=((\cF(E) \otimes \cE_2)\oplus \overline{\Delta_{\Theta}(\cF(E)\otimes
\cE_1)})
 \ominus \{ \Theta x \oplus \Delta_{\Theta}x : \;x \in \cF(E)\otimes \cE_1 \} .$$
Note that, if $\Theta$ is inner, we get $\Theta^*\Theta=U^*(q_1 \otimes
I_{\cE})U=I_{\cF(E)}\otimes q_1 $ and $\Delta_{\Theta}=0$. Thus, in this
case, $K(\Theta)=\cF(E)\otimes \cE_2$ and $H(\Theta)=(\cF(E)\otimes
\cE_2)\ominus \Theta(\cF(E)\otimes \cE_1)$.

We shall also write $P_{\Theta}$ for the projection from
$K(\Theta)$ onto $H(\Theta)$. We have the following.

\begin{theorem}\label{charep}
Let $\Theta$ be a characteristic operator function and let
$K(\Theta)$ and $H(\Theta)$ be as above. For every $a\in M$ and
$\xi \in E$ we define the operators $S_{\Theta}(\xi)$ and
$\psi_{\Theta}(a)$ on $\Delta_{\Theta}(\cF(E)\otimes \cE_1)$ by
$$S_{\Theta}(\xi)\Delta_{\Theta}g=\Delta_{\Theta}(T_{\xi}\otimes I_{\cE_1})g
,\;\;g\in \cF(E)\otimes \cE_1 $$
and
$$\psi_{\Theta}(a)\Delta_{\Theta}g=\Delta_{\Theta}(\varphi_{\infty}(a)\otimes
I_{\cE_1})g ,\;\;g\in \cF(E)\otimes \cE_1 .$$
Also, we define on $K(\Theta)$ the operators
$$V_{\Theta}(\xi)=(T_{\xi}\otimes I_{\cE_2})\oplus S_{\Theta}(\xi)
$$ and
$$\rho_{\Theta}(a)=(\varphi_{\infty}(a)\otimes I_{\cE_2})\oplus
\psi_{\Theta}(a).$$
Then
\begin{enumerate}
\item[(i)] $(S_{\Theta},\psi_{\Theta})$ and
$(V_{\Theta},\rho_{\Theta})$ define isometric covariant
representations of $E$ on $\overline{\Delta_{\Theta}(\cF(E)\otimes
\cE_1)}$ and $K(\Theta)$ respectively.
\item[(ii)] The space $K(\Theta)\ominus H(\Theta)$ is invariant
for $(V_{\Theta},\rho_{\Theta})$ and, thus, the compression of
$(V_{\Theta},\rho_{\Theta})$ to $H(\Theta)$, which we denote by
$(T_{\Theta},\sigma_{\Theta})$, is a completely contractive
covariant representation of $E$. Explicitely,
$$ T_{\Theta}(\xi)=P_{\Theta}V_{\Theta}(\xi)|H(\Theta)
\;,\;\;\xi\in E$$ and
$$\sigma_{\Theta}(a)=P_{\Theta}\rho_{\Theta}(a)|H(\Theta)\;,\;\;a\in
M.$$
\item[(iii)] The representation $(T_{\Theta},\sigma_{\Theta})$ is
completely non coisometric. It is a $C_{.0}$-representation if and
only if $\Theta$ is inner.
\end{enumerate}
\end{theorem}

The converse of the theorem above also holds; that is, every c.n.c
representation of $E$ gives rise to a characteristic operator
function. To construct it we now fix such a representation and
let $V,\pi, K, K_0,\mathcal{D}$ and $\mathcal{D}_*$ be as in
Section 2. Also, let $\rho_1$ be the restriction of $\pi$ to
$\mathcal{D}$ and $\rho_2$ be the restriction of $\pi$ (or
$\sigma$) to $\mathcal{D}_*$. (Again, we shall assume that
$\rho=\rho_1\oplus \rho_2$ is faithful. Otherwise a minor
technical correction is needed). We also write
$G=\mathcal{D}\oplus \mathcal{D}_*$.

Associated with the wandering subspaces $\mathcal{Q}_0$ and $\mathcal{D}$
 of $K$ we have the unitary operators
  $W(\mathcal{Q}_0)$ and $W(\mathcal{D})$ defined above (see the discussion
  preceeding Theorem~\ref{Theorem310MS98}). We write
  $Q_{\infty}$ for the projection of $K$ onto
  $L_{\infty}(\mathcal{Q}_0)$ (=the range of $W(\mathcal{Q}_0)$).

 Also recall that $u$ is an isometry from $\mathcal{Q}_0$ onto
 $\mathcal{D}_*$ that commutes with $\pi(a)$ for $a\in M$. It induces an
  isometry, written $I_{\cF(E)}\otimes u$
from $\cF(E) \otimes \mathcal{Q}_0$ onto $\cF(E)\otimes \cD_*$.

 We now write
\be\label{deftheta}
\Theta_T=(I_{\cF(E)}\otimes u) W(\mathcal{Q}_0)^*Q_{\infty}W(\cD)
:\cF(E) \otimes \cD
\rightarrow
\cF(E)\otimes \cD_* .
\ee

\begin{theorem}\label{charofT}
Given a c.n.c representation $(T,\sigma)$ of $E$ on a Hilbert
space $H$, the tuple
$(\Theta_T,\mathcal{D},\mathcal{D}_*,\rho_1,\rho_2 )$, defined
above, is a characteristic operator function in the sense of
Definition~\ref{charopfun}.
The representation is a $C_{.0}$-representation if and only if the
characteristic operator function is inner.

Moreover, the representation constructed from $\Theta_T$ as in
Theorem~\ref{charep} is unitarily equivalent to $(T,\sigma)$.
\end{theorem}

We also have the following.

\begin{theorem}\label{isom}
Suppose we start with a characteristic operator function \\
$(\Theta,\mathcal{E}_1,\mathcal{E}_2,\tau_1,\tau_2)$ and write
$(T_{\Theta},\sigma_{\Theta})$ for the c.n.c representation
constructed in Theorem~\ref{charep}. With this representation we
can associate the characteristic operator function
$(\Theta_T,\mathcal{D},\mathcal{D}_*,\rho_1,\rho_2)$ as in
Theorem~\ref{charofT}.

Then, the two characteristic operator functions are unitarily
equivalent in the sense that there are unitary operators
$W_1:\mathcal{E}_1 \rightarrow \mathcal{D}$ (intertwining $\tau_1$
and $\rho_1$) and $W_2:\mathcal{E}_2 \rightarrow \mathcal{D}_*$
(intertwining $\tau_2$ and $\rho_2$) such that
$$\Theta_T=(I_{\mathcal{F}(E)}\otimes W_2)\Theta (I_{\mathcal{F}(E)}\otimes
W_1^*). $$
\end{theorem}

The above results show that the characteristic operator functions
(or the elements $\widehat{\Theta}$ obtained from these) serve as
 complete invariants for c.n.c representations of $H^{\infty}(E)$.

We also state the following result which we know only for
$C_{.0}$-representations.
In this result we refer to compositions
$\Theta=\Theta_1\Theta_2$ where $\Theta$ is the inner
characteristic operator function associated with a
$C_{.0}$-representation and $\Theta_i$, $i=1,2$, is an inner characteristic
operator function but is not necessarily purely contractive. Two
such compositions $\Theta=\Theta_1\Theta_2=\Theta_1' \Theta_2'$
are said to be \emph{equivalent} if $\Theta_1'=\Theta_1(I\otimes
V_0)$ (and $\Theta_2'=(I\otimes V_0^*)\Theta_2$) for some unitary
operator $V_0$.

\begin{theorem}\label{invariant}
Let $(T,\sigma)$ a $C_{.0}$-representation of $E$ on $H$
(with $\sigma
\times T$ the associated representation of $H^{\infty}(E)$). Let
$\Theta$ be the inner characteristic operator function of this
representation.
Then there is a bijection between the subspaces of $H$ that are
$H^{\infty}(E)$-invariant and (equivalence classes of) factorizations
$\Theta=\Theta_1 \Theta_2$ of $\Theta$ as a composition of two
inner
characteristic operator functions (that are not necessarily purely contaractive).
\end{theorem}

\end{section}


\begin{thebibliography}{999}

\bibitem{AnD90} C. Anantharaman-Delaroche, \emph{On completely positive maps defined by an
irreducible correspondence}, Canad. Math. Bull. 33 (1990), 434-441.

\bibitem{Ar89} W.B. Arveson, \emph{Continuous analogues of Fock space},
Mem. Amer. Math. Soc. 80 (1989).

\bibitem{BDH88} M. Baillet, Y. Denizeau and J.-F. Havet, \emph{Indice
d'une esperance conditionelle}, Comp. Math. 66 (1988), 199-236.

\bibitem{BS00} B. V. R. Bhat and M. Skiede, \emph{Tensor product systems
of Hilbert modules and dilations of completely positive
semigroups}, Inf. Dim. Anal. Quantum Prob. and Rel. Topics 3
(2000), 519-575.

\bibitem{DP98} K. Davidson and D. Pitts, \emph{The algebraic structure
of non-commutative analytic Toeplitz algebras}, Math. Ann. 311
(1998), 275-303.

\bibitem{D74} R.G. Douglas, \emph{Canonical models}, in ``Topics in
operator theory", Math. Surveys 13 (Ed. C. Pearcy) (1974),
161-218.

\bibitem{KP03} D. Kribs and S. Power, \emph{Free semigroupoid algebras},
Preprint.

\bibitem{L94} E.C. Lance , \emph{Hilbert $C^{\ast}$-modules, A toolkit
for operator algebraists}, London Math. Soc. Lecture Notes series
210 (1995). Cambridge Univ. Press.

\bibitem{MM83} M. McAsey and P.S. Muhly, \emph{Representations of
non-self-adjoint crossed products}, Proc. London Math. Soc. 47
(1983), 128-144.

\bibitem{Mi89} J. Mingo, \emph{The correspondence associated to an inner
completely positive map}, Math. Ann. 284 (1989), 121-135.

\bibitem{MS98} P.S. Muhly and B. Solel, \emph{Tensor algebras over
$C^{\ast}$-correspondences (Representations, dilations and
$C^{\ast}$-envelopes)}, J. Funct. Anal. 158 (1998), 389-457.

\bibitem{MS99} P.S. Muhly and B. Solel , \emph{Tensor algebras, induced
representations, and the Wold decomposition}, Canad. J. Math. 51
(1999), 850-880.

\bibitem{MS00} P.S. Muhly and B. Solel, \emph{On the Morita equivalence
of tensor algebras}, Proc. London Math. Soc. 81 (2000), 361-377.

\bibitem{MS02} P.S. Muhly and B. Solel, \emph{Quantum Markov processes
(correspondences and dilations)}, Int. J. Math. 13 (2002), 863-906.

\bibitem{MS03} P.S. Muhly and B. Solel, \emph{The curvature and index of
completely positive maps}, Proc. London Math. Soc. 87 (2003),
748-778.

\bibitem{MSp03} P.S. Muhly and B. Solel, \emph{Hardy algebras,
$W^{\ast}$-correspondences and interpolation theory}, to appear in Math. Annalen.

\bibitem{MSCM} P.S. Muhly and B. Solel, \emph{On canonical models for
representations of Hardy algebras}. In preparation.

\bibitem{MSSp02} P.S. Muhly, M. Skeide and B. Solel, \emph{Representations of $\mathcal{B}^a(E)$, commutants of von Neumann bimodules, and product systems of Hilbert modules}. In preparation.

\bibitem{NF66} B. Sz-Nagy and C. Foias, \emph{Analyse Harmonique des
Operateurs de L'espace de Hilbert}, Akademiai Kiado (1966).

\bibitem{wP73} Wm. Paschke, \emph{Inner product modules of
B*-algebras}, Trans. Amer. Math. Soc. \textbf{182} (1973),
443-468.

\bibitem{Pet} J. Peters, \emph{Semi-crossed products of C*-algebras}, J.
Funct. Anal. 59 (1984), 498-534.

\bibitem{Pi97} M. Pimsner, \emph{A class of $C^{\ast}$-algebras
generalyzing both Cuntz-Krieger algebras and crossed products by
$\mathbb{Z}$}, in \textit{Free Probability Theory}, D. Voiculescu,
Ed., Fields Institute Comm. 12, 189-212, Amer. Math. Soc.,
Providence, 1997.

\bibitem{P86} S. Popa, \emph{Correspondences}, Preprint (1986).

\bibitem{Po89} G. Popescu, \emph{Characteristic functions for infinite
sequences of noncommuting operators}, J. Oper. Theory 22 (1989),
51-71.

\bibitem{Po91} G. Popescu, \emph{von Neumann inequality for
$B(\mathcal{H}^n)_1$}, Math. Scand. 68 (1991), 292-304.

\bibitem{R74} M.A. Rieffel, \emph{Morita equivalence for
$C^{\ast}$-algebras and $W^{\ast}$-algebras}, J. Pure Appl. Alg. 5
(1974), 51-96.

\bibitem{R74b} M.A. Rieffel, \emph{Induced representations of
$C^{\ast}$-algebras}, Adv. in Math. 13 (1974), 176-257.

\bibitem{mSp02} M. Skeide, \emph{Commutants of von Neumann modules, representations of $\cB^a(E)$  and other topics related to product systems of Hilbert modules}, to appear in the proceedings of the AMS Joint Summer Research Conference ``Advances in Quantum Dynamics", held at Mt. Holyoke College in the summer of 2002.

\bibitem{mSp03} M. Skeide, \emph{von Neumann modules, intertwiners and self-duality}, to appear in J. Operator Theory.

\end{thebibliography}
\end{document}